\documentclass[12pt,thmsa]{amsart}
\usepackage{amssymb, euscript}

\def\Cal{\mathcal}

\def\S{{\Cal S}}

\def\bbr{{\Bbb R}}

\def\bbn{{\Bbb N}}

\def\bbc{{\Bbb C}}

\def\supp{{\hbox{\rm supp}}}

\def\rn{\bbr^n}

\def\part{\partial}

\def\b{\beta}

\def\Gam{\Gamma}

\def\a{\alpha}

\def\Del{\Delta}

\def\vp{\varphi}

\def\sig{\sigma}
\def\lam{\lambda}

\def\t{\tau}

\newtheorem{theorem}{Theorem}[section]
\newtheorem{lemma}[theorem]{Lemma}

\newtheorem{corollary}[theorem]{Corollary}

\theoremstyle{remark}
\newtheorem{remark}[theorem]{Remark}

\numberwithin{equation}{section}

\newcommand{\be}{\begin{equation}}
\newcommand{\ee}{\end{equation}}

\newcommand{\bea}{\begin{eqnarray}}
\newcommand{\eea}{\end{eqnarray}}
\newcommand{\Bea}{\begin{eqnarray*}}
\newcommand{\Eea}{\end{eqnarray*}}

\begin{document}

\title[Spherical Means ]
{Spherical Means in Odd Dimensions\\ and EPD equations}

\author [B. Rubin]{Boris Rubin}
\address{
Department of Mathematics, Louisiana State University, Baton Rouge,
LA, 70803 USA}

\email{borisr@math.lsu.edu}

\thanks{The  research was supported in part by the  NSF grant DMS-0556157
and the Louisiana EPSCoR program, sponsored  by NSF and the Board of
Regents Support Fund.}





\begin{abstract} The paper contains a simple proof of the
Finch-Patch-Rakesh inversion formula for the spherical mean Radon
transform in odd dimensions. This transform arises in thermoacoustic
 tomography. Applications are given to the Cauchy problem for the Euler-Poisson-Darboux
 equation with initial data on the cylindrical surface. The argument relies on the idea
 of analytic continuation and known properties of Erd\'elyi-Kober fractional
 integrals.
\end{abstract}

\maketitle

\section{Introduction}

\setcounter{equation}{0}

We consider the spherical mean Radon transform  \be\label {mf}
(Mf)(\theta, t)=\frac{1}{\sig_{n-1}}\int_{S^{n-1}} f (\theta
-t\sig)\, d\sig,\ee $$ \theta \in S^{n-1}, \qquad t\in
\bbr_+=(0,\infty),$$ where $f$ is a smooth function supported inside
the unit ball $B=\{x \in \rn: |x|<1\}$, $S^{n-1}$ is the unit sphere
in $\rn$ with the area $\sig_{n-1}=2\pi^{n/2}/\Gam (n/2)$, and
$d\sig$ denotes integration against the usual Lebesgue measure on
$S^{n-1}$. The inversion
 problem  for this transform has attracted considerable attention in
 the last decade in view of new developments in thermoacoustic
 tomography; see \cite{AKK, AKQ, FPR, FHR, KK, Ku, PS1, PS2} and
 references therein. Explicit inversion formulas for $Mf$ in the ``closed form" are of
 particular interest. For $n$ odd, such formulas were obtained by Finch, Patch, and
 Rakesh  in \cite{FPR}. The corresponding formulas for  $n$ even were
 obtained by Finch,  Haltmeier, and  Rakesh in \cite{FHR}. An interesting inversion formula of different
 type,
 that covers both odd and even cases, was suggested by  Kunyansky \cite{Ku}; see also a survey paper
 \cite{KK}, where diverse inversion algorithms and related mathematical problems are discussed.

In spite of their elegance and ingenuity, most of explicit inversion
formulas for $Mf$ are still mysterious, and the basic ideas  behind
them are not completely understood. We observe, for example, that
the derivation for $n=3$  in \cite{FPR} relies on implementation of
delta functions and is not completely rigorous. The dimensions $n=5,
7, \ldots$ have been treated using a pretty complicated reduction to
the 3-dimensional case. The resulting inversion formula can be
written in the form \be\label{rf}  f(x)\!=\!c_n
\Del\int_{S^{n-1}}[D^{n-3}t^{n-2}\vp_\theta]\Big |_{t=|x\!-\!\theta
|}\, d\theta, \ee where $\vp_\theta (t)=(Mf)(\theta, t)$,
$$ \Del
\!=\!\sum\limits_{k=1}^n\frac{\partial^2}{\partial x_k^2}, \qquad
c_n\!=\!\frac{(-1)^{(n-1)/2}}{4\pi^{n/2-1}\,\Gam (n/2)},\qquad
D=\frac{1}{2t}\,\frac{d}{dt}\,;$$ cf. \cite[Theorem 3]{FPR}. Other formulas in
that theorem can be obtained in the framework of the same method;
see Remark \ref{rema}.

In the present article (Section 3) we suggest a simple and rigorous
proof of (\ref{rf}) that handles all odd $n$ simultaneously, so that
no reduction to $n=3$ is needed. The  idea is to treat the spherical
mean
 as a member of a certain analytic family of operators
associated to the Euler-Poisson-Darboux equation \be\label{epd}
\square_\a u \equiv \Del_x u - u_{tt} -\frac{n+2\a
-1}{t}\,u_{t}=0\ee  and invoke known facts from fractional calculus
\cite{SKM}
 about Erd\'elyi-Kober operators. Section 4 deals with applications.
 Here  we reconstruct the solution $u(x,t)$ of the equation $\square_\a
 u=\lam^2 u$, $ \lam \ge 0$, from the Cauchy
 data on the cylinder $S^{n-1} \times \bbr_+$. Section 5 contains
 comments and open questions.

No claim about the originality of the results presented in this
article is made, but it is felt that the elementary use of operators
of fractional integration to obtain them might appeal to the applied
mathematician. If the reader is not interested in applications in
Section 4, he may  skip subsections 2.2 and 2.4. These are not
needed for the basic Section 3.

 {\bf Acknowledgements.} I am grateful to professors Mark Agranovsky,  Peter
 Kuchment, and Larry Zalcman for very pleasant discussions and hospitality
 during my short visits to Bar Ilan and Texas A\&M Universities.

\section{Preliminaries}

\subsection{Erd\'elyi-Kober fractional integrals} We remind some
known facts; see, e.g., \cite[Sec. 18.1]{SKM}. For $Re \, \a >0$ and
$\eta \ge -1/2$, the Erd\'elyi-Kober fractional integral of a
function $\vp$ on $\bbr_+$ is defined by \be\label{eci}
(I^\a_\eta\vp)(t)=\frac{2t^{-2(\a+\eta)}}{\Gam (\a)}\int_0^t (t^2
-r^2)^{\a -1}r^{2\eta +1}\vp (r) \, dr,\qquad t>0.\ee For our
further needs, it suffices to assume that $\vp$ is infinitely smooth
and supported away from the origin. Then $I^\a_\eta\vp$ extends as
an entire function of $\a$ and $\eta$, so that
\be\label{ek0}I^0_\eta\vp=\vp,\ee \be\label{ek1}
(I^\a_\eta\vp)(t)=t^{-2(\a+\eta)}  D^m\,
t^{2(\a+m+\eta)}(I^{\a+m}_\eta\vp)(t), \qquad
D=\frac{1}{2t}\,\frac{d}{dt}\, ,\ee
  \be\label{ek2}
(I^{-m}_\eta\vp)(t)=t^{-2(\eta-m)}  D^m\, t^{2\eta}\vp(t),\ee where
$m$ is a nonnegative integer. The property \be D^m=t^{-1}\left
(\frac{d}{dt}\, \frac{1}{2t}\right )^m\, t\ee allows us to write
(\ref{ek1}) and (\ref{ek2}) in a different equivalent form. The
composition formula and the inverse operator are as follows:
\be\label{ek3} I^\b_{\eta+\a}I^\a_\eta \vp=I^{\a+\b}_\eta\vp, \qquad
(I^\a_\eta)^{-1}\vp=I^{-\a}_{\eta+\a}\vp.\ee

\subsection{The generalized Erd\'elyi-Kober fractional integrals}
Let $J_\nu$ and $I_\nu $ be  the Bessel function  and the modified
Bessel function of the first kind, respectively \cite{Er}. The
generalized Erd\'elyi-Kober operators
  are defined by
\be\label{defe}  J^\a_{\eta, \lam}\vp(t)=t^{-2(\a+\eta)}J^\a_{
\lam}t^{2\eta}\vp(t), \quad I^\a_{\eta,
\lam}\vp(t)=t^{-2(\a+\eta)}I^\a_{\lam}t^{2\eta}\vp(t),\ee where
$\lam \ge 0$, \bea\label{geci} \qquad J^\a_{\lam}\vp(t)\!&=&\!2^\a
\lam^{1-\a}\int_0^t (t^2\! -\!r^2)^{(\a -1)/2}\, J_{\a-1}(\lam
\sqrt{t^2 \!-\!r^2}\,) \vp (r) r\, dr,\\
\label{igeci}\qquad I^\a_{ \lam}\vp(t)\!&=&\!2^\a
\lam^{1-\a}\int_0^t (t^2\! -\!r^2)^{(\a -1)/2}\, I_{\a-1}(\lam
\sqrt{t^2 \!-\!r^2}\,) \vp (r) r\, dr; \eea see \cite{Lo1, Lo2},
\cite [Sec. 37.2]{SKM}. As above, we assume $\vp$ to be infinitely
smooth and supported away from the origin. Integrals (\ref{geci})
and (\ref{igeci}) are absolutely convergent if $Re\, \a>0$ and admit
analytic continuation to all complex $\a$ by the formulas
\bea\label{val} J^\a_{\lam}\vp\!&=&\!D^m
J_{\lam}^{\a+m}\vp\!=\!J^{\a+m} D^m\vp,
\\\label{vval} I^\a_{\lam}\vp\!&=&\!D^m
I^{\a+m}_{\lam}\vp\!=\!I_{\lam}^{\a+m} D^m\vp, \eea $ m \in \bbn$.
These follow from the well-known relation
$$\left ( \frac{1}{\t} \frac{d}{d\t}\right )^m \![\t^\nu J_\nu
(\t)]\!=\!\t^{\nu -m} J_{\nu -m} (\t)$$ (similarly for $I_\nu$).
Clearly, $J^\a_{\eta, 0} = I^\a_{\eta, 0}=I^\a_\eta$.  If $Re\, \a
>0$ and $ Re\, \b >0$, then, changing the
order of integration and using the formula 2.15.15(1) from
\cite{PBM},   we get
 \be I^\b_{\eta+\a, \lam}J^\a_{\eta,
\lam}\vp=I^{\a+\b}_\eta \vp.\ee The latter extends by analyticity to
all complex $\a$ and $\b$ and yields the inversion formula
\be\label{inve} (J^\a_{\eta, \lam})^{-1}f= I^{-\a}_{\eta+\a,
\lam}f.\ee By (\ref{vval}) and (\ref{defe}), this can also be
written as \be\label{laba} (J^\a_{\eta, \lam})^{-1}f=t^{-2\eta}D^m
I_{\lam}^{m-\a}t^{2(\eta+\a)}f=t^{-2\eta}D^m
t^{2(\eta+m)}I^{m-\a}_{\eta+\a, \lam}f.\ee

\subsection{The Euler-Poisson-Darboux
equation}

Consider the Cauchy problem for the Euler-Poisson-Darboux equation
(\ref{epd}):
 \be\label{epdk} \square_\a u=0, \quad u(x,0)=f(x),\quad
u_t(x,0)=0,\ee where $f$ belongs to the Schwartz space $\S(\rn)$;
see \cite{B1} for details. If $\a\ge (1-n)/2$, then (\ref{epdk}) has
a unique solution $ u(x,t)=(M_t^\a f)(x)$ where the operator $M_t^\a
$ is defined in the Fourier terms by $$ [M_t^\a f]^{\wedge}(y)=m_\a
(t|y|)\hat f (y),$$
$$ m_\a (\rho)=\Gam (\a+n/2) \,(\rho/2)^{1-\a-n/2}\,J_{n/2
+\a-1}(\rho).$$
 The
 operator $M_t^\a $ extends meromorphically to all complex $\a$ with the poles
 $ -n/2, -n/2-1, \ldots \,$. For $Re \, \a >0$, it is an integral
 operator \be\label{io}
 (M_t^\a f)(x)=\frac{\Gam (\a+n/2)}{\pi^{n/2}\Gam
 (\a)}\int_{|y|<1}(1-|y|^2)^{\a -1} f (x-ty)\, dy.\ee
In the case $\a=0$, $M_t^\a f$ is the spherical mean \be\label{mfs}
(M_t^0 f)(x)=\frac{1}{\sig_{n-1}}\int_{S^{n-1}} f (x -t\sig)\,
d\sig;\ee cf. (\ref{mf}). Passing to polar coordinates, one can
obviously represent $M_t^\a f$ as an Erd\'elyi-Kober  integral of
the spherical mean \be\label{rep} (M_t^\a f)(x)\!=\!\frac{\Gam
(\a\!+\!n/2)}{\Gam (n/2)}\, (I^\a_\eta \vp_x)(t), \qquad
\vp_x(t)\!=\!(M_t^0 f)(x),\ee with $\eta=n/2 -1$.

\subsection{The generalized Euler-Poisson-Darboux
equation}
 Consider the more general Cauchy
problem \be \label{lepdk}\square_\a u=\lam^2 u, \quad
u(x,0)=f(x),\quad u_t(x,0)=0, \ee where $f$ is a Schwartz function
and $\lam\ge 0$. If $\a\ge (1-n)/2$, then (\ref{lepdk}) has a unique
solution $ u(x,t)=(M_{t, \lam}^\a f)(x)$, where the operator $M_{t,
\lam}^\a $ is defined  as analytic continuation of the integral
\bea\label{lio}
 &&(M_{t, \lam}^\a f)(x)=\frac{(2/\lam)^{\a -1}\,\Gam (\a+n/2) }{\pi^{n/2}}\, t^{2-n-2\a}\\
 &\times&\int_{|y|<t} f (x-y)(t^2-|y|^2)^{(\a -1)/2}\, J_{\a  -1} (\lam \sqrt{t^2-|y|^2}\,)\,dy;\nonumber\eea
  see
\cite{B2} for details.  As above,  \be\label{mtl} (M_{t, \lam}^\a f)
(x)=\frac{\Gam (\a+n/2)}{\Gam
 (n/2)}\, (J^\a_{\eta,\lam }\vp_x)(t), \qquad \vp_x(t)\!=\!(M_t^0 f)(x),\ee
 $\eta=n/2 -1$, where $J^\a_{\eta,\lam }$ is
 the generalized
Erd\'elyi-Kober operator (\ref{defe}).

\subsection{More preparations} We restrict $M_t^\a f$ to
$x\!=\!\theta \!\in \!S^{n-1}$ and set \be\label{nef} (N^\a
f)(\theta, t)=t^{n+2\a -2}(M_t^\a f)(\theta).\ee In particular, for
$Re \, \a
>0$, owing to (\ref{io}), we have
 \be\label{naf} (N^\a f)(\theta, t)=\frac{\Gam (\a+n/2)}{\pi^{n/2}\Gam
 (\a)}\int_{\rn} f(y) \, (t^2-|y-\theta|^2)_+^{\a -1}\, dy\ee
where $(...)_+^{\a -1}$ has a standard meaning, namely, $(a-b)_+^{\a
-1}=(a-b)^{\a -1}$ if $a>b$ and $0$ otherwise. Given a function $F$
on the cylinder $S^{n-1} \times \bbr_+$, we denote \be\label{pf}
(PF)(x)= \frac{1}{\sig_{n-1}}\int_{S^{n-1}} F(\theta, |x-\theta |)\,
d\theta, \qquad x \in \rn,\ee which is a modification of the
back-projection (or dual) operator; cf. \cite {H,N}, where these
notions are used for the classical Radon transform. We  also invoke
Riemann-Liouville integrals \cite{SKM} \be\label{rli} (I_{-1}^\a
u)(s)=\frac{1}{\Gam (\a)}\int_{-1}^s (s-t)^{\a -1} u(t)\, dt, \qquad
u\in C^\infty [-1,1].\ee
 The integral (\ref{rli})  is absolutely convergent when $Re \, \a >0$ and
 extends by analyticity to all $\a\in \bbc$, so that \be\label{eba}(I_{-1}^{-m}
 u)(s)=(d/ds)^m \,u(s), \qquad m=0,1,2, \ldots \,.\ee

\begin{lemma} Let $B=\{x \in \rn: |x|<1\}$. For $Re \, \a >0$ and any integrable function $f$
supported in $B$ we have \be\label {anal} (PN^\a
f)(x)=c_\a\int_B\frac{f(y)}{|x-y|^{1-\a}} \,(I_{-1}^\a u)(h)\,
dy,\qquad x\in B, \ee where
$$ c_\a\!=\!\frac{\Gam (n/2)\,\Gam (\a+n/2)}{2^{1-\a}\,\pi^{(n+1)/2}\,\Gam
  ((n\!-\!1)/2)}, \quad  h\!=\!\frac{|x|^2\!-\!|y|^2}{2|x-y|}, \quad u(t)\!=\!(1-t^2)^{(n-3)/2}.$$
\end{lemma}
\begin{proof} By (\ref{naf}) and (\ref{pf}), changing the order of
integration, we get \bea (PN^\a f)(x)\!&=&\!\frac{\Gam
(\a+n/2)}{\sig_{n-1}\,\pi^{n/2}\,\Gam
 (\a)}\int_{S^{n-1}}\!\!\!d\theta\!\int_B \!f(y) (|x\!-\!\theta |^2\!-\!|y\!-\!\theta
 |^2)_+^{\a -1}\, dy\nonumber\\ &=&\!\frac{\Gam
(\a+n/2)}{\sig_{n-1}\,\pi^{n/2}\, \Gam
 (\a)}\int_B f(y)\,k_\a (x,y)\, dy,
\nonumber\eea where \bea k_\a (x,y)&=&\int_{S^{n-1}}
(|x|^2-|y|^2-2\theta \cdot (x-y))_+^{\a -1}\, d\theta\nonumber
\\&=&\frac{\sig_{n-2}}{(2|x-y|)^{1-\a}}\int_{-1}^h (h-t)^{\a -1} (1-t^2)^{(n-3)/2}\,
dt.\nonumber\eea This gives the result.\end{proof}

\section{Inversion of the Spherical Mean for $n$ Odd} Let $C^\infty (B)$ be the space
of $C^\infty$-functions on $\rn$ supported in $B$; $f \in C^\infty
(B)$. By (\ref{nef}),
 (\ref{rep}), and (\ref{ek2}), analytic continuation of $N^\a f$ at
 $\a=3-n$ has the form $$ (N^{3-n} f)(\theta, t)=\frac{\Gam
 (3-n/2)}{\Gam (n/2)}\, D^{n-3}\, t^{n-2}\vp_\theta (t), \qquad \vp_\theta
 (t)=(Mf)(\theta, t).$$
If $n=2k+3, \; k=0,1, \ldots \,$, then (\ref{eba}) yields
$$(I_{-1}^{3-n} u)(h)=(d/dh)^{2k} (1-h^2)^k=(-1)^k k!= (-1)^{(n-3)/2}\Gam
(n-2).$$ Hence, analytic continuation of (\ref{anal}) at $\a=3-n$ is
\be\label{prp} P[D^{n-3}t^{n-2}\vp_\theta](x)\!=\!c \,(I^2f)(x),
\qquad c\!=\!2(-1)^{(n-3)/2}\Gam^2 (n/2)/\pi,\ee where $x \in B$ and
\be\label{rpo}(I^2f)(x)=\frac{\Gam (n/2 -1)}{4\pi^{n/2}}\int_B
\frac{f(y)\, dy}{|x-y|^{n-2}}\ee is the Riesz potential of order $2$
\cite{SKM}. The latter can be inverted by the Laplacian, and simple
calculations yield \be\label{my} f(x)\!=\!c_n
\Del\int_{S^{n-1}}[D^{n-3}t^{n-2}\vp_\theta]\Big |_{t=|x\!-\!\theta
|}\, d\theta, \quad c_n\!=\!\frac{(-1)^{(n-1)/2}}{4\pi^{n/2-1}\,\Gam
(n/2)}.\ee

\begin{remark}\label{rema} Formula (\ref{my}) coincides (up to notation) with the third
formula in \cite[Theorem 3]{FPR}. Other formulas in that theorem can
be similarly obtained from (\ref{anal}) if the latter is applied to
$\Del f$ instead of $f$. Here we take into account that $I^2\Del
f=-f$, because $\supp f$ is separated from the boundary of $B$.
\end{remark}

\section{An inverse problem for the EPD equation}

Let $\a\ge (1-n)/2, \; \lam \ge 0$. Suppose we know the trace
$u(\theta,t)$ of the solution of the Cauchy problem
\be\label{cpo}\square_\a u=\lam^2 u, \quad u(x,0)=f(x),\quad
u_t(x,0)=0,\ee on the cylinder $\{(\theta,t): \theta \in S^{n-1}, \,
t\in \bbr_+\}$ and want to reconstruct the initial function $f$ in
the space $C^\infty (B)$. This can be easily done using (\ref{prp})
and the Erd\'elyi-Kober operators. Indeed, by (\ref{mtl}),
$$u(\theta,t)\!=\!(M_{t, \lam}^\a f)(\theta)\!=\!\frac{\Gam (\a\!+\!n/2)}{\Gam (n/2)}\,
(J^\a_{\eta,\lam } \vp_\theta)(t), $$ where
$\vp_\theta(t)\!=\!(M_t^0 f)(\theta)$, $ \eta\!=\!n/2\! -\!1$. Then
by (\ref{inve}), for $ u_\theta (t) \equiv u(\theta,t)$ we have
$$
\vp_\theta\!=\!\frac{\Gam (n/2)}{\Gam
(\a\!+\!n/2)}\,(J^\a_{\eta,\lam })^{-1}u_\theta\!=\!\frac{\Gam
(n/2)}{\Gam (\a\!+\!n/2)}\, I^{-\a}_{\eta+\a, \lam}u_\theta.$$ Now,
 (\ref{prp}) yields
$$\frac{\Gam (n/2)}{\Gam (\a\!+\!n/2)}\,
 P[D_t^{n-3}t^{n-2}I^{-\a}_{\eta+\a, \lam}u_\theta]\!=\!c \,I^2f, \quad c\!=\!2(-1)^{(n-3)/2}\Gam^2 (n/2)/\pi, $$
and therefore, \be f=\frac{(-1)^{(n-1)/2}\, \pi}{2\Gam(n/2)\,\Gam
(\a\!+\!n/2)}\,\Del P[D_t^{n-3}t^{n-2}I^{-\a}_{\eta+\a,
\lam}u_\theta].\ee The operator $I^{-\a}_{\eta+\a, \lam}$ can be
replaced by any expression from (\ref{laba}).

It remains to note that once $f$ is known, the solution $u(x,t)$ of
the equation  $\square_\a u=\lam^2 u$ can be reconstructed from the
trace $u(\theta,t)$  by the formula $u(x,t)=(M_{t, \lam}^\a f)(x)$;
see (\ref{lio}).

\section{Comments}

{\bf 1.} It is a challenging open problem to appropriately adjust
our method to the case when $n$ is even and give alternative proof
of the corresponding inversion formulas from \cite{FHR} and
\cite{Ku}.

{\bf 2.} Formula (\ref{prp})  provokes the following

\noindent{\bf Conjecture. }  {\it Let $n\ge 3$ be odd. A function
$\vp_\theta (t) \!\equiv\!\vp (\theta,t)$ belongs to the range of
the operator $f \to (Mf)(\theta,t)$, $f\! \in\! C^\infty (B)$, if
and only if $P[D_t^{n-3}t^{n-2}\vp_\theta]$ belongs to the range
$I^2 [ C^\infty (B)]$ of the  potential (\ref{rpo}).}

The ``only if" part follows immediately from  (\ref{prp}). The ``if"
part requires studying injectivity of the back-projection operator
$P$, which is of independent interest; cf. \cite{Ru2}, where
injectivity and inversion of the dual Radon transform is studied in
the general context of affine Grassmann manifolds. Various
descriptions of the range of the spherical mean transform and many
related results can be found in \cite{AKQ, AK, FR}.

{\bf 3.} It is worth noting that for $n=3$, the inversion formula
for $Mf$ becomes elementary if  $f$ is a radial function, i.e.
$f(x)\equiv f_0(|x|)$.
\begin{lemma} If $f\in L^1_{loc} (\rn), \; f(x)\equiv f_0(|x|)$,
then $(Mf)(\theta, t)\equiv F_0 (t)$, where \be\label{chv} F_0
(t)=\frac{2^{n-3}\, \Gam (n/2)}{\pi^{1/2}\, \Gam
((n-1)/2)}\int_{|1-t|}^{1+t} f_0 (r) \,[a(r,t)]^{n-3}\, r\, dr\ee
 $a(r,t)=[r^2 -(1-t)^2]^{1/2}[(1+t)^2-r^2]^{1/2}/4$ being the area
of the triangle with sides $1,t,r$.
\end{lemma}
\begin{proof}
\bea (Mf)(\theta, t)&=&\frac{1}{\sig_{n-1}}\int_{S^{n-1}}f_0
(|\theta -t\sig|)\, d\sig\nonumber
\\&=&\frac{\sig_{n-2}}{\sig_{n-1}}\int_{-1}^1
f_0(\sqrt{1+t^2-2ts})\,(1-s^2)^{(n-3)/2}\, ds \nonumber\eea and
(\ref{chv}) follows.
\end{proof}
\begin{corollary} If $n=3$ and $f$ is supported in $B$, then (\ref{chv}) yields
\be (a) \quad F_0 (t)=\frac{1}{2t} \int_{1-t}^1 f_0 (r) r\, dr \quad
\text{\rm if}\quad 0<t\le 1,\ee and \be (b) \quad F_0
(t)=\frac{1}{2t} \int_{t-1}^1 f_0 (r) r\, dr \quad \text{\rm
if}\quad 1\le t<2.\ee In the case (a), \be\label{fo1} f_0
(r)=\frac{2}{r}\left [\frac{d}{dt} (tF_0 (t)) \right ]_{t=1-r}. \ee
 In the case (b), \be\label{fo2} f_0
(r)=-\frac{2}{r}\left [\frac{d}{dt} (tF_0 (t)) \right ]_{t=1+r}. \ee
\end{corollary}
\begin{remark} We note that in (\ref{fo1}) and  (\ref{fo2})  it
is not necessary to know $F_0 (t)$ {\it for all} $t \in (0,2)$ as in
(\ref{my}). It suffices to know it only for $t \in (0,1)$ or $t \in
(1,2)$. It would be interesting to obtain similar inversion formulas
for $(Mf)(\theta, t)$ in the general case, when $(Mf)(\theta, t)$ is
known only for $(\theta ,t)  \in S^{n-1}\times (0,1)$ or $(\theta
,t)  \in S^{n-1}\times  (1,2)$, as in the radial case.

\end{remark}

\end{document}